\documentclass[12pt]{article}
\usepackage{amsthm,amsfonts,amsmath,amssymb}
\usepackage{graphicx}
\usepackage{ tabls}
\usepackage{marvosym}

\usepackage{pstricks,pst-node,pst-tree, pstcol}

\newcommand{\ts}{{\tt s}}
\newcommand{\tm}{{\tt h}}

\begin{document}
\title{The Unlucky Door} 
\author {Sasha Gnedin\thanks{\tt A.V.Gnedin@uu.nl }}

\maketitle

\section{The game}

Monte and Conie play the famous Three-Door Game. The quiz-team hides the prize behind one of the doors.
Conie, who does not know where the prize was hidden, is asked to choose one of the doors as a first guess. 
Monte, who saw where the prize was hidden, will reveal then
one of the doors which does not conceal the prize,
but never the door chosen by Conie. Finally, Conie will be offered to either hold her choice or to switch to another yet unrevealed door.
Conie wins if her final choice is the door which hides the prize.

See the book by Rosenhouse for history and variations of the problem \cite{Rosenhouse}.
The Three-Door-Game in proper sense of the game theory, as interaction of two actors, appeared before in   \cite{Gill}, \cite{GameTh},
\cite{Haggstrom}.
In this note we shall focus on a combinatorial aspect of the game and 
a possibility of cooperative play for a certain design with four doors.

Let us label the doors $1,2,3$ and think of four moves in the game. 
The first move is simple: the quiz-team hides the prize behind door $p$. 
On the second move Conie chooses door $x$. On the third move Monte offers a switch to door $y\neq x$ by revealing a door which is not $p$.
On the fourth move Conie chooses $z$ from $x$ and $y$: if she decides to hold her initial choice then $z=x$ and if she decides to switch
then $z=y$. She wins if $z=p$.

The strategy of the quiz-team is just the action of hiding the prize. We shall think of this move as a move of nature.

What can Conie do? On move two she chooses $x$, and on move four makes her finial choice  $z$  which depends on both $x$ and $y$.
For instance, she may first guess $x=1$ and then decide to hold if $y=2$ and to switch if $y=3$.
Her strategy can be labeled like e.g. $2\tm\ts$ which means the algorithm 
``first choose door 2, then hold if a switch to door with smaller number is offered
and switch if a switch to a door with larger number is offered''. There are twelve such strategies $1\ts\ts, 1\tm\tm$ etc. 
It is a good exercise to write down them all!

What can Monty do? He knows $p$, and on the second move observes $x$.  If $p\neq x$ then his decision is forced 
to offer a switch to another door,
and if $p=x$ he can choose out of two options: to offer a switch to a door with smaller or larger number. 
We can label his six strategies by sequences  like $212$ which indicate the switch offer $y$ as reaction 
on match of the first Conies guess with door $p=1,2,3$, respectively. Monte has six strategies in total.

It is important to understand that when Conie and Monte fix the way they will play the game, that is choose their strategies, the 
course of the game is completely determined by $p$. We can consider then the course of the game 
as the work of computer program which has input parameter $p$. In particular, for given profile of two 
strategies, one  of Monte and one of Conie,
the value of $p$ 
determines if Conie wins or not.

\section{Win or lose}

\begin{verse}
{\it One, two, three\\
That's how elementary it's gonna be...}
\end{verse}

\noindent
Conie has  strategy $1\ts\ts$ which wins for $p=2,3$ and loses for $p=1$, so wins in two cases out of three,
no matter how Monte plays.
More generally, each always-switching strategy $x\ts\ts$ loses in case $p=x$ and wins in two other cases.

It is intuitively obvious that no  Conie's strategy can win in all three  cases $p=1,2,3$ for any given Monte's play.
One explanation for that is the following
\vskip0.5cm
\noindent
{\bf The Unlucky Door Theorem.} {\it For every strategy $S$ of Conie there exists at least one  door $u=u(S)$, which depends neither 
on   $p$ nor on Monte's strategy, such that
Conie loses when $p=u$ for every strategy of Monte.}
\vskip0.3cm

The proof is simple. If $S=1\ts\ts$, then the prize is not won for $p=1$ and we can take $u=1$. If $S=1\tm\tm$ then Conie holds
$x=1$ whichever happens and so we can take $u=2$ or $3$.
If $S=1\ts\tm$ then Conie will not switch to $y=3$, but for $p=3\neq 1=x$ Monte will offer precisely $y=3$, 
so we can take $u=3$. Similarly, $u=2$ for  $S=1\ts\tm$.
The general principle to find the unlucky door is: 
there is always a door $u$ which is never Conie's choice at the forth 
move, whichever strategy of Monte.

\vskip0.3cm
We note in passing, that the number of winning cases may depend on Monte's play.
For instance
for strategy $1\tm\ts$~  Conie wins for $p=1,3$ versus any strategy $2\cdots$, and wins only for $p=3$ versus Monte's $3\cdots$.
The complete matrix of winning cases is found in \cite{GameTh},\cite{Haggstrom}.

Using unlucky doors we can readily show that every strategy of Conie is weakly dominated by some always-switching strategy.
For strategy $S$ one just takes $u\ts\ts$. Since $u\ts\ts$ loses only for $p=u$, and $S$ loses then too,
each time $S$ wins the strategy $u\ts\ts$ wins as well.

\vskip0.3cm

We see that whichever Conie does, Monte cannot play a strategy to make Conie sure winner, 
for all locations of the prize.
In particular, if Conie plays some always-switching strategy, e.g. $1\ts\ts$, she wins in 2 cases which is 
the maximum possible.

It is a well-known and obvious (and at the same time counter-intuitive)  fact that,
if the quiz-team places the prize at random by rolling a
three-sided symmetric die then every always-switching strategy wins with probability $2/3$, for every strategy of Monte.
The existence of unlucky door implies that no other strategy of Conie can have under this randomization
higher winning probability, because with probability at least $1/3$ the prize will be behind the unlucky door.

\section{Four door designs}

\begin{verse}
{\it
The hard part is learning about it\\
The hard part is breaking through to the truth\\
The hard part is learning to doubt what you read\\
What you hear, what you see on the news.}
\end{verse}

\noindent
We wish to extend the game to the case  with four doors and one prize.
 Suppose the first two moves are as in the basic version: the prize is hidden behind door $p\in\{1,2,3,4\}$ by the quiz-team, 
and door $x\in \{1,2,3,4\}$ is chosen by Conie.
On the third move, which we will consider with two plausible designs,
Monty reveals two doors as not hiding the prize and offers a switch to door $y$, which is not $x$.
On the fourth move Conie can either hold $x$ or switch to $y$, winning if her final choice $z$ is $p$.

\vskip0.5cm
\noindent
{\bf Revealing two doors at once.} We start with a simpler design. After Conie chooses $x$, the light is switched off and on.
Conie sees two doors revealed as not hiding the prize, and she is left with the  dilemma to switch to $y$ or to hold $x$.

Conie has in this game more strategies like 1\ts\tm\ts, meaning ``choose door $x=1$ then hold it for $y=3$,
and switch from it to $y$ for $y=2,4$''.  However, precisely the same argument we used for the basic game with three doors shows that for every
her strategy $S$ there exists at least one unlucky door $u=u(S)$ such that $u$ is never her final choice. When $p=u$ the game is 
lost whichever Monte's play.

On the other hand, every always-switching strategy like e.g. $1\ts\ts\ts$~ wins in three cases out of four.
Therefore this is the maximal possible number of cases $p\in\{1,2,3,4\}$ won. 
If a four-sided symmetric die is rolled to 
randomize $p$, Conie wins with any always-switching strategy with probability $3/4$, and this is the highest possible.

Exactly as in the case of three doors, Monte cannot help to win for sure. 
The same is true for extensions of the game to arbitrary number of doors $n\geq 3$, provided Monte reveals $n-2$ doors at once.

\vskip0.5cm
\noindent
{\bf Revealing in sequence. } A minor change in the design of the third move can change the situation dramatically.
If two out of four doors are revealed in a succession observed by Conie, 
 Monte can signal Conie complete information about the location of the prize. So if the game of 
Monte and Conie is cooperative, they can agree about a common strategy to make Conie sure winner.

The sets of strategies of both players is now large enough, as Monty can reveal 2 doors in sequence, and 
Conie can take final action  depending on both $x$ and the sequence.

We show one possible signaling algorithm by constructing a function 
$$f:(p,x)\mapsto (r_1,r_2)$$ 
which 
assigns to each value of $(p,x)$ two sequentially revealed door numbers $(r_1,r_2)$ which are neither $p$ nor $x$.

Take $x=1$, and consider a correspondence $(p,x)\mapsto (r_1,r_2)$ 
given by 
$$1\mapsto 23,~~ 2\mapsto 34,~~3\mapsto 42,~~4\mapsto 32.$$
Think of the doors 1,2,3,4 as arranged in the cyclic order. For the values $p=1,2,3$ we take for $r_1,r_2$ two admissible, 
next to $p$ clockwise 
 doors, and for $p=4$ two next counter-clock wise doors.
The function is extended to the values $x=2,3,4$ by the clockwise rotation of the labels.

One corollary of this construction is that for $n=4$ there is no analogue of {\it universally} unlucky door for each
Conie's strategy, such that 
behind the door  the prize is never found whichever $p$ and strategy of Monte.

Now suppose Monty uses the algorithm encoded in the function $f$ to reveal two doors, one-by-one. 
If the rule is comminicated to Conie or she has derived it from her previous experience,
she can, using the rule and the door $x$ she initially chose, determine inambiguously the location
of prize $p$, thus win the prize in all cases.

If $f$ with the described bijective property varies from one round to another in a way known to Conie, 
she will be still sure winner. The parameters $(p,x)$ can be in this case arbitrarily jointly randomized, in dependent or independent way,
without changing the sure outcome.
Randomizing $x$ and agreeing in advance on $f$, perhaps with changing $f$ from round to round,
the actors need not communicate during the game to impress e.g. the quiz-team with a long series of wins.

If $p$ is uniformly random  and Monte picks, one-by-one, a door to reveal at random from all admissible options, then
the strategy SLM (switch at the last minute) of choosing, say, $x=1$ and always switching as two unrevealed doors remain
yields the  winning probability $3/4$, which in this scenario is the maximal possible  \cite{Rosenhouse}.

\end{document}